\documentclass[11pt]{article}
\usepackage{amsmath}
\usepackage{amssymb}
\usepackage{amsthm}
\usepackage{geometry}
\usepackage{amscd}
\newtheorem{thm}{Theorem}[section]
\newtheorem{prop}[thm]{Proposition} %%Delete [thm] to re-start numbering
\newtheorem{cor}[thm]{Corollary} %%Delete [thm] to re-start numbering

\newtheorem{lem}[thm]{Lemma}

\theoremstyle{definition}
\newtheorem{defn}[thm]{Definition}
\newtheorem{rem}[thm]{Remark}

\AtBeginDocument{\renewcommand{\refname}{Bibliography}}
\begin{document}

\begin{center}
\textbf{\Large{Generalized Inner-Outer Factorizations in non commutative Hardy Algebras}}
\end{center}

\begin{center}
Leonid Helmer
\end{center}

%\tableofcontents
%\addcontentsline{toc}{section}{Abstract}\begin{abstract}
\begin{abstract}
Let $H^{\infty}(E)$ be a non commutative Hardy algebra, associated with a $W^*$-correspondence $E$.
In this paper we construct factorizations of inner-outer type of the elements of $H^{\infty}(E)$ represented via the induced representation, and of the elements of its commutant. These factorizations generalize the classical inner-outer factorization of elements of $H^\infty(\mathbb{D})$. Our results also generalize some results that were obtained by several authors in some special cases.
\end{abstract}

\section{Introduction}

In this work we describe the general version of the inner-outer factorization in non commutative Hardy algebras. Recall that the Hardy algebra $H^{\infty}(\mathbb{D})$ is identified with the algebra $H^{\infty}:=H^{\infty}(\mathbb{T}):=L^{\infty}(\mathbb{T})\cap H^{2}(\mathbb{T})$, where $H^2=H^2(\mathbb{T})$ is the Hardy Hilbert space, and we consider $H^{\infty}$ as the algebra of multiplication operators $M_{\phi}$ acting on the Hilbert space $H^2$ by $f\mapsto \phi f$. Then the function $\Theta\in H^{\infty}$ is called inner if the operator $M_{\Theta}$ is isometric and the function $g\in H^{\infty}$ is called outer if the operator $M_{g}$ has a dense range in $H^2$. The classical theorem says that every $f\in H^{\infty}$ admits a unique inner-outer factorization $f=f_if_o$, where $f_i$ is an inner function, called also the inner part of $f$, and $f_o$ is an outer, called the outer part of $f$. Analogous factorizations hold in the Hardy spaces $H^p$, $p\geq 1$. In particular, every $f\in H^2$ admits an inner-outer factorization o$f=\Theta g$, with $f_i\in H^{\infty}$ and $f_o\in H^2$.
Further, any $z$-invariant subspace of the form $\mathcal{M}_f=\vee\{z^nf: n=0,1,...\}$ has the representation $\mathcal{M}_f=f_{i}H^2$. The classical Beurling' theorem says that every $z$-invariant subspace $\mathcal{M}$ has a representation $\mathcal{M}=\Theta H^2$ for a suitable inner function $\Theta$, \cite{Beu}.
A full treatment of the classical theory both from the function theoretic and the operator theoretic point of view can be found in \cite{Nik} and \cite{NF}.

Before we introduce the non commutative Hardy algebras note that the classical algebra $H^{\infty}(\mathbb{D})$ can be viewed as the ultraweak closure of the operator algebra generated by the unilateral shift on the Hilbert space $l^2=l^2(\mathbb{Z}_+)$. In \cite{Po89} this was generalized by G. Popescu to the ultraweakly closed non commutative operator algebras generated by $d$ shifts, $d\geq 1$, denoted $\mathcal{F}^\infty$. In \cite{APo}, \cite{Po2}, \cite{Po3}, \cite{Po4} Arias and Popescu developed the theory of inner-outer factorization in $\mathcal{F}^\infty$. In \cite{DavP2} Davidson and Pitts developed analogous theory, with some differences, in the context of the free semigroup algebra $\mathcal{L}_d$, which, in fact, coincides with $\mathcal{F}^\infty$.
In \cite{KriPow} Kribs and Power considered the case of free semigroupoids algebras $\mathcal{L}_G$, and developed the theory of inner-outer factorization in these algebras.

In this work we develop our version of the inner-outer factorization in non commutative Hardy algebras $H^{\infty}(E)$ associated with a given $W^*$-correspondence $E$. These algebras were introduced in 2004 by P. Muhly and B. Solel in ~\cite{MuS3} (see also ~\cite{MuS2}), and generalize the classical Hardy algebra $H^\infty$, the algebra $\mathcal{F}^\infty$ of Popescu, free semigroups algebras, free semigroupoids algebras and some others.

Let $E$ be a $W^*$-correspondence over a $W^*$-algebra $M$, (\cite{ManT}, \cite{Pa}), that is a right Hilbert $W^*$-module $E$ over $M$, which is made into a $M$-$M$ - bimodule  by some $*$-homomorphism $\phi: M\rightarrow \mathcal{L}(M)$, where $\mathcal{L}(M)$ is the algebra of all the adjointable operators on $E$. This $W^*$-correspondence defines another $W^*$-correspondence $\mathcal{F}(E)$ over the same algebra $M$, which is defined to be the direct sum $M\oplus E\oplus E^{\otimes 2}\oplus...$ of the internal tensor powers of $E$. $\mathcal{F}(E)$ is called the full Fock space and in fact is a $W^*$-correspondence with the left action of $M$ denoted by $\phi_{\infty}$, which is a natural extension of $\phi$ to a representation of $M$ in the algebra of adjointable operators on $\mathcal{F}(E)$. Note that the space $\mathcal{L}(E)$, for any $W^*$-correspondence $E$, is a $W^*$-algebra. The non commutative Hardy algebra of a correspondence $E$ is by definition the weak$^*$-closure in $\mathcal{L}(\mathcal{F}(E))$ of the algebra spanned by the operators of the form $T_{\xi}$, $\xi\in E$, where $T_\xi(\eta):=\xi\otimes \eta$ and $\phi_{\infty}(a)$, $a\in M$. In fact Muhly and Solel defined this Hardy algebra as the weak$^*$ closure of the noncommutative tensor algebra $\mathcal{T}_+(E)$. The algebra $\mathcal{T}_+(E)$ was defined first in ~\cite{MuS2} as the norm closed (nonselfadjoint) algebra spanned by the same set of generators, and it generalizes the noncommutative disc algebra $\mathcal{A}_n$ of Popescu, which in its turn is a noncommutative generalization of the classical disc algebra.

In this work we view the algebra $H^{\infty}(E)$ as acting on a Hilbert space via an induced representation $\rho$ and write it $\rho(H^{\infty}(E))$. Thus we consider the questions of inner-outer factorization for the case of the Hardy algebra $\rho(H^{\infty}(E))$. A key tool that we will need and use here is the general version of Wold decomposition proved first in ~\cite{MuS1}. We start with the inner-outer factorization of a vector of the underlying Hilbert space, and then we obtain the inner-outer factorization of an element of the commutant of the algebra $\rho(H^{\infty}(E))$. Here we use that fact that, as in the abstract theory of shifts, an inner operator is a partial isometry in the commutant of the algebra, generated by a shift. Further, we translate the Beurling theorem of Muhly and Solel in ~\cite{MuS1} to our language. It follows from the concept of duality for $W^*$-correspondences, developed in  ~\cite{MuS3}, every algebra $\rho(H^\infty(E))$ can be thought of as the commutant of the Hardy algebra of another correpondence, called the dual of $E$. Using this concept we construct factorization of an element of $\rho(H^\infty(E))$ which holds in our setup.

\section{Preliminaries and Setting}

We start by recalling the notion of a $W^*$-correspondence. For a general theory of Hilbert $C^*$- and $W^*$-modules  we use \cite{Lan}, \cite{ManT} and the original paper \cite{Pa}. Here we only note that by a Hilbert $W^*$-module we always mean a self dual module over a $W^*$-algebra (see \cite[Ch. 3]{ManT}).

Let $\phi: M\rightarrow \mathcal{L}(E)$ be a normal $*$-homomorphism. In what follows we always assume that $\phi$ is unital. Then we obtain on $E$ the structure of a bimodule over $M$. We shall call it a $W^*$-correspondence over the $W^*$-algebra $M$. More generally, let $N$ and $M$ be $W^*$-algebras and let $E$ be a Hilbert $M$-module. Assume that we are given a left action of $N$ on $E$, that is, we are given normal $*$-homomorphism
$\phi:N\rightarrow \mathcal{L}(E)$. This homomorphism can be regarded as a ``generalized homomorphism" from $N$ to $M$. Such an $N$-$M$-bimodule $E$ will be called a correspondence from $N$ to $M$. Every $W^*$-correspondence $E$ has the structure of a dual Banach space \cite{Pa}. This topology is usually called the $\sigma$-topology, \cite{MuS3}.

Every Hilbert space $H$, where the inner product is taken to be linear in the second variable, is a $W^*$-module and a $W^*$-correspondence over $\mathbb{C}$ in a
natural way.

Let $E$ and $F$ be $W^*$-correspondences over $W^*$-algebras $M$ and $N$ respectively. The left action of $M$ on $E$ will be denoted as usual by $\phi$ and the left action of $N$ on $F$ by $\psi$, thus, $\psi:N\rightarrow \mathcal{L}(F)$ is a normal $*$-homomorphism.
\begin{defn}\label{Isomorphism of W^*-corresp} An isomorphism of $E$ and $F$ is a pair $(\sigma,\Phi)$ where

1) $\sigma:M\rightarrow N$ is an isomorphism of $W^*$-algebras;

2) $\Phi:E\rightarrow F$ is a vector space isomorphism preserving the $\sigma$-topology, and which is also

~~(a) a bimodule map, $\Phi(\phi(a)xb)=\psi(\sigma(a))\Phi(x)\sigma(b)$, $x\in E$ $a,b\in M$, and

~~(b) $\Phi$ ``preserves" the inner product, $\langle\Phi(x),\Phi(y)\rangle=\sigma(\langle x,y\rangle)$, $x,y\in E$.
\end{defn}

Let $E$ be a $W^*$-correspondence over a $W^*$-algebra $M$ with a left action defined as usual by a normal $*$-homomorphism $\phi$.
For each $n\geq 0$, let $E^{\otimes n}$ be the self-dual internal tensor power (balanced over $\phi$, \cite{MuS3}). So, $E^{\otimes n}$ itself turns out to be a $W^*$-correspondence in a natural way, with the left action $\xi \mapsto \phi_n(a)\xi=(\phi(a)\xi _{1})\otimes ...\otimes \xi_n$, $\xi=\xi_1\otimes...\xi_n\in E^{\otimes n}$, and with an $M$-valued inner product as in the internal tensor product construction. For example, on $E^{\otimes 2}=E\otimes _{\phi} E$, we define
$$\langle \xi_1\otimes \xi_2,\eta_1\otimes \eta_2\rangle=
\langle \xi_2,\phi(\langle\xi_1,\eta_1\rangle)\eta_2\rangle.$$

We form the full Fock space $\mathcal{F}(E)=\sum^{\oplus}_{n\geq 0}E^{\otimes n}$, where $E^{\otimes 0}=M$ and the direct sum taken in the ultraweak sense (see \cite{Pa}). This is a $W^*$-correspondence with left action given by
$\phi _{\infty}:M\rightarrow \mathcal{L}(\mathcal{F}(E))$, where $\phi _{\infty}(a)=\sum _{n\geq 0} \phi_n(a)$. The $M$-valued inner product on $\mathcal{F}(E)$ is defined in an obvious way.

For each $\xi\in E$ and each $\eta\in \mathcal{F}(E)$, let $T_{\xi}: \eta \mapsto \xi \otimes \eta$ be a creation operator on
$\mathcal{F}(E)$. Clearly, $T_{\xi}\in \mathcal{L}(\mathcal{F}(E))$.

\begin{defn}\label{Defn of tensor and Hardy alg} Given a $W^*$-correspondence $E$ over a $W^*$-algebra $M$.

1) The norm closed subalgebra in $\mathcal{L}(\mathcal{F}(E))$, generated by all creation operators $T_{\xi}$, $\xi\in E$, and all operators $\phi _{\infty}(a)$,  $a\in M$, is called the tensor algebra of $E$. It is denoted by $\mathcal{T}_{+}(E)$.

2) The Hardy algebra $H^{\infty}(E)$ is the ultra-weak closure
of $\mathcal{T}_{+}(E)$.
\end{defn}

When $M=E=\mathbb{C}$ then $\mathcal{F}(E)=l^2(\mathbb{Z}_+)$. The algebra
$\mathcal{T}_{+}(\mathbb{C})$ is the algebra of analytic Toeplitz operators with continuous symbols, so it can be identified with the disc algebra $A(\mathbb{D})$. The algebra
$H^{\infty}(\mathbb{C})$, in this case, is $H^{\infty}(\mathbb{D})$.
If $M=\mathbb{C}$ and we take $E=\mathbb{C}^n$, an $n$-dimensional Hilbert space, then $\mathcal{T}_+(\mathbb{C}^n)$ is the non commutative disc algebra $\mathcal{A}_n$, studied by Popescu and others, and $H^{\infty}(\mathbb{C}^n)=\mathcal{F}^{\infty}_n$, the Hardy algebra of Popescu. This algebra can be identified with the free semigroup algebra $\mathcal{L}_n$ studied by Davidson and Pitts.

Let $\pi:M\rightarrow B(H)$ be a normal representation of a $W^*$-algebra $M$ on a Hilbert space $H$ and let $E$ be a $W^*$-correspondence over $M$. As it can be easy verified, the $W^*$-internal tensor product $E\otimes_{\pi}H$ is a Hilbert space. The representation $\pi^{E}:\mathcal{L}(E)\rightarrow B(E\otimes _{\pi}H)$ defined by

$$\pi^{E}: S\mapsto S\otimes I_H, \,\,\,\ \forall S\in \mathcal{L}(E).$$
is called the induced representation (in the sense of Rieffel). If $\pi$ is a faithful normal representation then $\pi^E$ maps $\mathcal{L}(E)$ into
$B(E\otimes _{\pi}H)$ homeomorphically with respect to the ultraweak topologies, \cite[Lemma 2.1]{MuS3}.

The image of $H^{\infty}(E)$ under an induced representation is defined as follows. Let $\pi:M\rightarrow B(H)$ be a faithful normal representation. For a $W^*$-correspondence $E$ over $M$ let
$\pi^{\mathcal{F}(E)}$ be the induced representation of $\mathcal{L}(\mathcal{F}(E))$ in $B(\mathcal{F}(E)\otimes_{\pi}H)$.
Then the induced representation of the Hardy algebra $H^{\infty}(E)$ is the restriction
\begin{equation}\label{Induced repr of Hardy alg}
\rho:=\pi^{\mathcal{F}(E)}|_{H^{\infty}(E)}:H^{\infty}(E)\rightarrow B(\mathcal{F}(E)\otimes_{\pi}H).
\end{equation}
This restriction is an ultraweakly continuous representation
of ${H^{\infty}(E)}$ and the image $\rho(H^{\infty}(E))$ is an ultraweakly closed subalgebra of $B(\mathcal{F}(E)\otimes_{\pi}H)$. We shall refer to $\rho$ as the representation induced by $\pi$. Later, when we discuss several representation of $H^{\infty}(E)$ that are induced by different representations $\pi$, $\sigma$ etc. of $M$, we shall write $\rho_{\pi}$, $\rho_\sigma$ etc.

So, $\rho(H^{\infty}(E))$ acts on
$\mathcal{F}(E)\otimes_{\pi}H$ and $\rho$ is defined by
$$\rho:X\mapsto X\otimes I_H,\,\,\,\forall{X}\in H^{\infty}(E).$$

Note that the notion of the induced representation generalizes the notion of pure isometry (i.e. an isometry without a unitary part) in the theory of a single operator.

We will frequently use the following result of Rieffel ~\cite[Theorem 6.23]{Rie}.
The formulation here is in a form convenient for us (~\cite[p. 853]{MuS1}).
\begin{thm}\label{Rieffel thm}. Let $E$ be a $W^*$-correspondence over the algebra $M$ and $\pi:M\rightarrow B(H)$ be a normal faithful representation of $M$ on the Hilbert space $H$. Then the operator $R$ in $B(E\otimes_{\pi}H)$ commutes with $\pi^{E}(\mathcal{L}(E))$ if and only if $R$ is of the form $I_E\otimes X$, where $X\in \pi(M)'$, i.e., $\pi^{E}(\mathcal{L}(E))'=I_E\otimes \pi(M)'$.
\end{thm}

\subsection{Covariant representations.}

\begin{defn}
Let $E$ be a $W^*$-correspondence over a $W^*$-algebra $M$.

(1) By a covariant representation of $E$, or of the pair $(E,M)$, on a Hilbert space $H$, we mean
a pair $(T, \sigma)$, where $\sigma:M\rightarrow B(H)$ is a nondegenerate normal $*$-homomorphism, and
$T$ is a bimodule (with respect to $\sigma$) map $T:E\rightarrow B(H)$, that is a linear map such
that $T(\xi a)=T(\xi)\sigma(a)$ and $T(\phi(a)\xi)=\sigma(a)T(\xi)$, $\xi\in E$ and $a\in M$. We require also that $T$ will be continuous with respect to the $\sigma$-topology on $E$ and the ultraweak topology on $B(H)$.

(2) The representation $(T,\sigma)$ is called (completely) bounded, (completely) contractive, if so is the map $T$. For a completely contractive covariant representation we write also c.c.c.r.

\end{defn}

The operator space structure on $E$ to which this definition refers is the one which comes from the embedding of $E$ into its so-called linking algebra $\mathfrak{L}(E)$, see ~\cite{MuS2}.

In this work we will consider only isometric covariant representations. A covariant representation $(V,\sigma)$ is said to be isometric if $V(\xi)^*V(\eta)=\sigma(\langle \xi,\eta\rangle)$, for every $\xi,\eta\in E$. Every isometric covariant representation $(V,\pi)$ of $E$ is completely contractive, see ~\cite[ Corollary 2.13]{MuS2}.

As an important example let $\rho=\pi^{\mathcal{F}(E)}|_{H^{\infty}(E)}$ be an induced representation of the Hardy algebra $H^{\infty}(E)$. For the representation $\sigma$ set
$$\sigma=\pi^{\mathcal{F}(E)}\circ\phi_{\infty},$$
and set
$$V(\xi)=\pi^{\mathcal{F}(E)}(T_{\xi}),\,\,\,\ \xi\in E.$$
\begin{defn}\label{Def of ind cov rep}
The pair $(V,\sigma)$ is called the covariant representation induced by $\pi$, or simply the induced covariant representation (associated with $\rho$).
\end{defn}

It is easy to check that $(V,\sigma)$ in the above definition is isometric, hence, is completely contractive.

Let $(T, \sigma)$ be a c.c.c.r. of $(E,M)$ on the Hilbert space $H$ as above.
With each such representation we associate the operator $\tilde{T}:E\otimes _{\sigma}H\rightarrow H$, that on the elementary tensors is defined by
$$\tilde{T}(\xi \otimes h):=T(\xi)(h).$$
$\tilde{T}$ is well defined since $T(\xi a)=T(\xi)\sigma(a)$.
In ~\cite{MuS2} Muhly and Solel show that the properties of $\tilde{T}$ reflect the properties of the covariant representation $(T,\sigma)$.
They proved that $(\alpha)$ $ \tilde{T}$ is bounded iff $T$ is
completely bounded, and in this case $\|T\|_{cb}=\|\tilde{T}\|$; $(\beta)$ $\tilde{T}$ is
contractive iff $T$ is completely contractive; and $(\gamma)$ $\tilde{T}$ is an isometry  iff
$(T,\sigma)$ is an isometric representation.
A simple calculation gives us the intertwining relation
\begin{equation}\label{intetwining for T}
\tilde{T}\sigma^{E}\circ\phi(a)=\tilde{T}(\phi(a)\otimes I_H)=\sigma (a)\tilde{T}, \,\,\ \forall
a\in A.
\end{equation}

In the following theorem we collect two basic facts concerning the theory of representations of $W^*$-correspondences and of their tensor algebras.
\begin{thm}\label{bijection (T,sigma-T_tildeAndIntegratedFormOfRepr} (\cite[Lemma 2.5 and Theorem 2.9]{MuS3}) Let $E$ be any $W^*$-correspondence over an algebra $M$.  Then

1) There is a bijective correspondence $(T,\sigma)\leftrightarrow \tilde{T}$ between all
c.c.c.r. $(T,\sigma)$ of $E$ on a Hilbert space $H$ and contractive operators $\tilde{T}:E\otimes_{\sigma}H\rightarrow H$ that satisfy the relation (\ref{intetwining for T}). Let $\tilde{T}:E\otimes_{\sigma}H\rightarrow H$ be a contraction that satisfies the relation (\ref{intetwining for T}). Then the associated covariant representation is the pair $(T,\sigma)$, where $T$ is defined by
$T(\xi)h:=\tilde{T}(\xi\otimes h)$, $h\in H$ and $\xi\in E$.

2) Let $E$ be a $W^*$-correspondence over the algebra $M$ and let
$(T,\sigma)$ be a c.c.c.r. of $(E,M)$ on a Hilbert space $H$. Then for every such representation there exists a completely contractive representation $\rho:\mathcal{T}_+(E)\rightarrow B(H)$ such that
$\rho(T_{\xi})=T(\xi)$ for every $\xi\in E$ and $\rho(\phi_{\infty}(a))=\sigma(a)$ for every $a\in M$. Moreover, the correspondence $(T,\sigma)\leftrightarrow\rho$ is a bijection between the set of all c.c.c.r. of $E$ and all completely contractive representations of $\mathcal{T}_+(E)$ whose restrictions to $\phi_{\infty}(M)$ are continuous with respect to the ultraweak topology on $\mathcal{L}(\mathcal{F}(E))$.
\end{thm}

Restricting our attention to isometric covariant representation, we have the following.

\begin{lem} \label{descript_of_isometric_repr}(\cite[Lemma 2.1]{MuS1}.)
Let $(V,\sigma)$ be any isometric covariant representation  of the $W^*$-correspondence $E$ on a Hilbert space $H$. Then the associated
isometry $\tilde{V}:E\otimes_{\sigma}H\rightarrow H$ is an isometry that satisfy the relation
$\tilde{V}\sigma^{E}\circ\phi(a)=\sigma (a)\tilde{V}$, $\forall a\in M$, and with range equal to the closed linear span of $\{V(\xi)h:\xi\in E,a\in M\}$. Conversely, given an isometry $\tilde{V}:E\otimes_{\sigma}H\rightarrow H$ that satisfies the above
intertwining relation, then the associated covariant representation is the pair $(V,\sigma)$, where $V$ is defined by
$V(\xi)h:=\tilde{V}(\xi\otimes h)$, $h\in H$ and $\xi\in E$.
\end{lem}

The representation $\rho$ of $\mathcal{T}_+(E)$ that corresponds to the covariant representation $(T,\sigma)$ is called the integrated form of $(T,\sigma)$ and denoted by $\sigma\times T$. In its turn, the representation $(T,\sigma)$ is called the desintegrated form of $\rho$.
Preceding results show that, given a normal representation $\sigma$ of
$M$, the set of all completely contractive representations of $\mathcal{T}_+(E)$ whose restrictions to $\phi_{\infty}(M)$ is given by $\sigma$ can be
parameterized by the contractions $\tilde{T}\in B(E\otimes_{\sigma}H, H)$, that satisfy the relation $(\ref{intetwining for T})$.

In this notations, the induced representation $\rho_\pi$ is an integrated form of the $(V,\sigma)$, the covariant induced representation of $E$ from Definition $\ref{Def of ind cov rep}$.

In \cite{MuS3} it was shown that, if
the representation $(T,\sigma)$ of $(E,M)$ is such that
$\|\tilde{T}\|<1$, then the integrated form $\sigma\times T$
extends from $\mathcal{T}_+(E)$ to an ultraweakly continuous representation of $H^{\infty}(E)$. For a general $(T,\sigma)$, the question when such an extention is possible is more delicate, see about this \cite{MuS8}.

Let $(V,\sigma)$ be an isometric covariant representation of a general $W^*$-correspondence $E$ on a Hilbert space $G$. For every $n\geq 1$ write $(V^{\otimes n},\sigma)$ for the isometric covariant representation of $E^{\otimes n}$ on the same space $G$ defined by the formula $V^{\otimes n}(\xi_1\otimes...\otimes\xi_n)=V(\xi_1)\cdots V(\xi_n)$, $n\geq 1$. The associated isometric operator $\tilde{V}_n:E^{\otimes n}\otimes_\sigma G\rightarrow G$ (which is called the generalized power of $\tilde{V}$), satisfies the identity
$\tilde{V}_n\sigma^{E^{\otimes n}}\circ \phi_n= \tilde{V}_n (\phi_n\otimes I_{G})=\sigma \tilde{V}_n$. In this notation $\tilde{V}=\tilde{V}_1$.

For each $k\geq 0$ write $G_k$ for $\bigvee\{V(\xi_1)\cdots V(\xi_k)g:\xi_i\in E,g\in G\}$ (with $G_0=G$). Clearly, $G_k=\overline{\tilde{V}_k(E^{\otimes k}\otimes_{\sigma}G_0)}$. Write $R_k$ the projection of $G_0$ onto $G_k$ and let $P_k=R_{k}-R_{k+1}$ and $R_{\infty}=\wedge_kR_k$. Thus,
$R_k=\sum_{l\geq k}P_l+ R_{\infty}$ is a projection of $G_0$ onto $G_k$ and $R_0=I_{G_0}$.
According to \cite{MuS1}, the formula $L(x)=\tilde {V}_1(I_1\otimes x)\tilde {V}_1^*$ defines a normal endomorphism of the commutant $\sigma(M)'$ and its n-th iterate is $L^n(x)=\tilde {V}_n(I_n\otimes x)\tilde {V}_n^*$. Here $I_n=I_{E^{\otimes n}}$.
Simple calculation shows that $L^n(P_m)=P_{n+m}$ and $R_n=\tilde {V}_n\tilde {V}_n^*=L^n(I)$.

An isometric covariant representation $(V,\sigma)$ is called fully coisometric if $R_1=L(I_{G_0})=I_{G_0}$. Muhly and Solel proved the following Wold decomposition theorem (\cite[Theorem 2.9]{MuS1}):

\begin{thm}Let $(V,\sigma)$ be an isometric covariant representation of $W^*$-correspondence $E$ on a Hilbert space $G_0$. Then $(V,\sigma)$ decomposes into a direct sum $(V_1,\sigma_1)\oplus (V_2,\sigma_2)$ on $G_0=H_1\oplus H_2$, where $(V_1,\sigma_1)=(V,\sigma)|_{H_1}$ is an induced representation and $(V_2,\sigma_2)=(V,\sigma)|_{H_2}$ is fully coisometric. Further, this decomposition is unique in the sense that if $K\subseteq G_0$ reduces $(V,\sigma)$ and the restriction $(V,\sigma)|_K$ is induced (resp. fully coisometric) then $K\subseteq H_1$ (resp. $K\subseteq H_2$).
\end{thm}
From this theorem it follows immediately that $(V,\sigma)$ is an induced representation if and only if $R_{\infty}=\wedge_kR_k=0$.

With $(V,\sigma)$ we may associate the ``shift" $\mathfrak{L}$, that acts on the lattice of $\sigma (M)$-invariant subspaces of $G$, and is defined as a geometric counterpart of the endomorphism $L$. In a more details, let $\mathcal{M}\in Lat(\sigma (M))$, then we set

\begin{equation}\label{DefnOfGeneralizedShift}
\mathfrak{L}(\mathcal{M}):=\bigvee \{V(\xi)k : \xi \in E, k\in \mathcal{M}\}.
\end{equation}
The $s$-power $\mathfrak{L}^s(\mathcal{M})$ is defined in the obvious way (with $\mathfrak{L}^0(\mathcal{M})=\mathcal{M})$).

The subspace $\mathcal{M}\in Lat(\sigma(M))$, as well as its projection $P_{\mathcal{M}}\in \sigma(M)'$, is called wandering with respect to $(V,\sigma)$, if the subspaces  $\mathfrak{L}^s(\mathcal{M})$, $s=0,1,...$, are mutually orthogonal.
Write $\sigma'$ for the restriction $\sigma|_{\mathcal{M}}$, where $\mathcal{M}$ is wandering. Then the Hilbert space $E^{\otimes s}\otimes_{\sigma'}\mathcal{M}$ is isometrically isomorphic (under the generalised power $\tilde{V}_s$) to $\mathfrak{L}^s(\mathcal{M})$. Hence, we obtain an isometric isomorphism
$$\mathcal{F}(E)\otimes_{\sigma'}\mathcal{M}\cong \sum ^{\oplus}_{s\geq 0}\mathfrak{L}^s(\mathcal{M}).$$

In these notations we have
$G_k=\mathfrak{L}^k(G_0)\cong E^{\otimes k}\otimes_{\sigma}G_0\cong\sum^{\oplus}_{l\geq k}E^{\otimes l}\otimes_{\pi}H$, with $H$ as the wandering subspace and $\sigma'=\pi$.

\section{Generalized inner-outer factorization}

In this section we describe a general version of the theory of inner-outer factorization for an arbitrary element $g\in \mathcal{F}(E)\otimes_\pi H$ and arbitrary elements of commutant $\rho(H^{\infty}(E))'$, and then we deduce some natural version of factorization of elements of $\rho(H^{\infty}(E))$, where $\rho=\rho_\pi$ denotes the representation of $H^{\infty}(E)$ on $\mathcal{F}(E)\otimes_\pi H$, induced by the faithful normal representation $\pi$. Although most of our constructions are correct in the general case we assume in the following that the space $H$ of the representation $\pi$ is separable (see Remark $\ref{Cardinality of ampliation}$).

Before we start let $S$ be a unilateral shift acting on the Hilbert space $H$ and let $\mathcal{M}\subset H$ be an $S$-invariant subspace. Write $\mathcal{M}_0:=\mathcal{M}\ominus S(\mathcal{M})$ and $H_0:=H\ominus S(H)$ for the wandering subspaces of $S|_{\mathcal{M}}$ and of $S$ correspondingly. Then one of the main points in the proofs of the classical theorems of Beurling, Halmos and Lax on invariant subspaces of $S$ is that $\dim \mathcal{M}_0\leq \dim H_0$.

In our situation let us consider $G=\mathcal{F}(E)\otimes_\pi H$ as the left $H^{\infty}(E)$-module with the action defined by $X\cdot g:=\rho_\pi(X)g$, for any $X\in H^{\infty}(E)$ and $g\in G$.
Thus, in this language every $\rho_\pi(H^\infty(E))$-invariant subspace $\mathcal{M}\subseteq G$ defines an $H^\infty(E)$-submodule in $G$. Note that in this case $End(G)$ - the set of all the endomorphisms of this module is nothing but $\rho_\pi(H^\infty(E))'$. The covariant representation $(V,\sigma)$, associated with the induced representation $\rho_\pi$, defines the generalized shift $\mathfrak{L}$. Hence, we need to compare the wandering subspaces $G\ominus \mathcal{L}(G)$ and  $\mathcal{M}_0:=\mathcal{M}\ominus \mathcal{L}(\mathcal{M})$. More precisely, we need to compare the representations of $M$ on $H$ and on the $\mathcal{M}_0$. This is done in the following proposition

\begin{prop}\cite[Proposition 4.1]{MuS1}
Let $\mathcal{M}$ be a $\rho_\pi(H^\infty(E))$-invariant subspace of $\mathcal{F}(E)\otimes_\pi H$ and let $(V,\sigma)$ be the associated covariant representation of $(E, M)$. If $\mathcal{M}_0=\mathcal{M}\ominus \mathfrak{L}(M)$ is the $\mathfrak{L}$- wandering subspace in $\mathcal{F}(E)\otimes_\pi H$, then the restriction $\sigma|_{\mathcal{M}_0}$ is unitarily equivalent to a subrepresentation of $\pi$ if and only if there is a partial isometry in $\rho_\pi(H^\infty(E))'$ with final space $\mathcal{M}$
\end{prop}

As the induced covariant representation $(V,\sigma )$ is a natural generalization of a pure isometry, that is of a shift operator, a partial isometry in $\rho_\pi(H^\infty(E))'$ was called an inner operator, \cite{MuS1}. In our work we shall generalize this definition, and shall use this term for a suitable isometric operator which intertwines representations of $M$.

In fact we use the modules's language only to emphasize the analogy with the classical theory of shifts. Instead of this we shall constantly use the language of generalized shift $\mathfrak{L}$, associated with the induced covariant representation $(V,\sigma)$.

%%%%%%%%%%%%%%%%%%%%%%%%%%%%%%%%%%%%%%%%%%%%%%%%%%%%%%%%%%%%%%%%%%%%%%%%%%%%%%%%%%%%%%%%%%%%%%%%%%%%%%%%5
%%%%%%%%%%%%%%%%%%%%%%%%%%%%%%%%%%%%%%%%%%%%%%%%%%%%%%%%%%%%%%%%%%%%%%%%%%%%%%%%%%%%%%%%%%%%%%%%%%%%%%%%
\subsection{Inner-outer factorization of elements of $\mathcal{F}(E)\otimes_\pi H$}
We turn to the inner-outer factorization of a vector in $\mathcal{F}(E)\otimes_\pi H$. To this end we prove a Beurling type theorem for a cyclic $(V,\sigma)$-invariant subspace generated by this vector, i.e. for subspaces of the form $\mathcal{M}_g=\overline{\rho(H^{\infty}(E))g}$, where $g\in G_0:=\mathcal{F}(E)\otimes_\pi H$ is arbitrary.

Write $P_g$ for the projection $P_{\mathcal{M}_g}$ onto $\mathcal{M}_g$.
Clearly, $\mathcal{M}_g$ is a $\rho(H^{\infty}(E))$-invariant subspace, $P_g\in \sigma(M)$ and the restriction $(V,\sigma)|_{\mathcal{M}_g}$ is an induced isometric covariant representation, as follows from \cite[Proposition 2.11]{MuS1}.
In particular, $\mathcal{M}_g\in Lat(\sigma(M))$ and the subspace $\mathfrak{L}(\mathcal{M}_g)$ is well defined.
Set

\begin{equation}
\mathcal{N}_g:= \mathcal{M}_g\ominus \mathfrak{L}(\mathcal{M}_g).
\end{equation}

By $Q_g$ we denote the orthogonal projection of $G_0$ on $\mathcal{N}_g$. Then
$Q_g=P_g-L(P_g)$ is the wandering projection associated with the restricted representation $(V,\sigma)|\mathcal{M}_g$. Since $L^k(Q_g)\perp L^s(Q_g)$, $k\neq s$, and since $(V,\sigma)|\mathcal{M}_g$ is induced we obtain the decomposition
$\mathcal{M}_g=\sum^{\oplus}_{k\geq 0} L^k(Q_g)\mathcal{M}_g$. Equivalently, $\mathfrak{L}^k(\mathcal{N}_g)\perp \mathfrak{L}^s(\mathcal{N}_g)$, $k\neq s$, and $$\mathcal{M}_g=\mathcal{N}_g\oplus \mathfrak{L}(\mathcal{N}_g)\oplus\cdots .$$

We set $g_0:=Q_{g}g\in \mathcal{N}_g$.

\begin{lem}\label{wandering supsp for M_g}
$\mathcal{N}_g=\overline{\sigma (M)g_0}$.
\end{lem}

\noindent{\textbf{Proof.}}

Since $\sigma(a)g_0\in \mathcal{N}_g$ for each $a\in M$, then
$\overline{\sigma (M)g_0}\subset \mathcal{N}_g$. Let $z\in \mathcal{N}_g\ominus \overline{\sigma(M)g_0}$.
Then $z\perp \sigma(a)g_0$ for each $a\in M$
and in particular $z\perp g_0$.
Write $g=g_0+(g-g_0)$. Since $g-g_0\in \mathfrak{L}(\mathcal{M}_g)$, we get
$z\perp g-g_0$. For each $k\geq 1$, $V^{\otimes k}(\xi)g\in \mathfrak{L}^k
(\mathcal{M}_g)\subset \mathcal{M}_g\ominus \mathcal{N}_g$. So,
$z\perp V^{\otimes k}(\xi)g$ for each $k\geq 1$, $\xi\in E^{\otimes k}$.
It follows that $z\perp \mathcal{M}_g$ and then $z=0$, $\forall z\in \mathcal{N}_g\ominus \overline{\sigma(M)g_0}$, i.e.
$\mathcal{N}_g= \overline{\sigma(M)g_0}$. $\square$

\begin{rem}\label{two form of cyclic subsp.}
It is easy to see that every element of the form $\xi\otimes h\in E^{\otimes k}\otimes_\pi H$, for every $k\geq 1$ and every $\xi\in E^{\otimes k}$, is a wandering vector.
\end{rem}

For each $a\in M$ we set

\begin{equation}
\tau (a)=\langle \sigma (a)g_0,g_0 \rangle =\langle (\phi _{\infty}(a)\otimes I) g_0,g_0\rangle .
\end{equation}

This defines a positive ultraweakly continuous linear functional $\tau$ on $M$.
Since $\pi$ is assumed to be faithful, we can view $\tau$ as
defined on $\pi (M)\subset B(H)$.

Hence, there is a sequence $\{h_i\}\subset H$ with $\sum _i \|h_i\|^2\leq \infty$,
such that

\begin{equation}
\tau (a)=\sum _i\langle \pi (a)h_i,h_i\rangle=\langle\sigma(a)g_0,g_0\rangle .
\end{equation}

This sequence $\{h_i\}$ can be viewed as an element of the space $H^{(\infty)}=H\oplus H\oplus ...$ and we write
$h_{\tau}=\{h_i\}$ for it to indicate that it is defined by the functional $\tau$. For each $a\in M$ we define an operator
(ampliation of $\pi$) $\hat{\pi}(a)=diag (\pi (a))\in B(H^{(\infty)})$, acting
on $H^{(\infty)}$ by: $\hat{\pi}(a)k=\{\pi (a)k_i\}$ where $k=\{k_i\}\in H^{(\infty)}$.
Then, for $h_{\tau}$ we have
$$\tau(a)=\langle \hat{\pi} (a)h_{\tau},h_{\tau}\rangle = \sum _i\langle \pi (a)h_i,h_i\rangle =
\langle \sigma(a)g_0,g_0 \rangle .$$

Set
\begin{equation}\label{K_tau}
K_{\tau}:=\overline{\hat{\pi}(M)h_{\tau}}\subseteq H^{(\infty)},
\end{equation}
and define the operator $w_0 :H^{(\infty)}\rightarrow \mathcal{N}_g$, by

\begin{equation}\label{def of w_0}
\hat{\pi} (a)h_\tau \mapsto (\phi_{\infty}(a)\otimes I)g_0=\sigma(a)g_0,
\end{equation}
and $w_0=0$ on $H^{(\infty)}\ominus K_\tau$. Since $\langle \hat{\pi} (a)h_{\tau},h_{\tau}\rangle=\langle \sigma(a)g_0,g_0 \rangle$, $w_0$ is a well defined partial isometry from $H^{(\infty)}$ onto $\mathcal{N}_g$.
Taking $a=1\in M$, we get $w_0h_{\tau}=g_0$, and we see that
$$w_0(\hat{\pi}(a)h_\tau)=\sigma(a)g_0=\sigma(a)w_0(h_{\tau}).$$

Since the sets $\{\hat{\pi}(a)h_\tau:a\in M\}$ and $\{\sigma(a)g_0:a\in M\}$ are dense in $K_\tau$ and $\mathcal{N}_g$ respectively we obtain that $w_0$
intertwines $\hat{\pi}$ and $\sigma$:
\begin{equation}\label{intertwining for w_0}
w_0\hat{\pi}(a)=\sigma(a)w_0,\ \forall{a}\in M.
\end{equation}
We conclude:
\begin{prop}\label{Properties of w_0} The operator $w_0$ is a partial isometry intertwining $\hat{\pi}$
and $\sigma $, with initial subspace
$K_{\tau}$, and with $\mathcal{N}_g$ as final
subspace.
\end{prop}

Now let us consider the space $G_0^{(\infty)}=G_0\oplus G_0\oplus...$ and identify it with the space
$\mathcal{F}(E)\otimes_{\hat{\pi}} H^{(\infty)}$.
As usual we identify $M\otimes_{\hat{\pi}}H^{(\infty)}$ with $H^{(\infty)}$ and set
$\hat{\mathfrak{L}}(H^{(\infty)}):=\mathfrak{L}(H)\oplus \mathfrak{L}(H)\oplus ...
=\mathfrak{L}(H)^{(\infty)}$.
Hence, for $G_0^{(\infty)}$ we can write the decomposition

\begin{equation}
G_0^{(\infty)}=
H^{(\infty)}\oplus \hat{\mathfrak{L}}(H^{(\infty)})
\oplus \hat{\mathfrak{L}}^2(H^{(\infty)})\oplus ...
\end{equation}

Write $\hat{\rho}$ for the induced representation $\rho_{\hat{\pi}}$ of $H^{\infty}(E)$ on $\mathcal{F}(E)\otimes_{\hat{\pi}}H^{(\infty)}$, $X\mapsto X\otimes I_{H^{(\infty)}}$, $X\in H^{\infty}(E)$. Thus, $\hat{\rho}$ is an ampliation of the induced representation $\rho$ on $\mathcal{F}(E)\otimes_\pi H$. Then the associated isometric covariant representation of $E$ on
$\mathcal{F}(E)\otimes_{\hat{\pi}}H^{(\infty)}$ is the pair $(\hat{V},\hat{\sigma})$ where
$\hat{V}(\xi)=T_{\xi}\otimes I_{H^{(\infty)}}$ and
$\hat{\sigma}(a)=\phi_{\infty}(a)\otimes I_{H^{(\infty)}}$.
We call it an ampliation of $(V,\sigma)$.
Similarly we define the covariant representations $(\hat{V}^{\otimes k},\hat{\sigma})$ and the associated operators $(\hat{V}_k)^{\sim}$, $k\geq 1$.

For each $k\geq 0$ we identify $\hat{\mathfrak{L}}^k(K_{\tau})$ with $E^{\otimes k}\otimes_{\hat{\pi}}K_\tau$ and $\mathfrak{L}^k(\mathcal{N}_g)$ with $E^{\otimes k}\otimes_{\sigma}\mathcal{N}_g$. So, $\sum_{k\geq 0}\hat{\mathfrak{L}}^k(K_{\tau})=\sum_{k\geq 0}E^{\otimes k}\otimes_{\hat{\pi}}K_{\tau}=\mathcal{F}(E)\otimes_{\hat{\pi}}K_{\tau}$ and $\mathcal{M}_g=\sum_{k\geq 0}E^{\otimes k}\otimes_{\sigma}\mathcal{N}_g=\mathcal{F}(E)\otimes_{\sigma}\mathcal{N}_g$

Write $\sigma'$ for the restriction $\sigma|_{\mathcal{N}_g}$. Consider the restriction $\rho_\pi(H^{\infty}(E))|_{\mathcal{M}_g}$ and
the induced representation $\sigma'^{\mathcal{F}(E)}(H^{\infty}(E))$ of $H^{\infty}(E)$ on $\mathcal{M}_g$. Then for every
$k\geq 0$ and $\xi\otimes z\in E^{\otimes k}\otimes_{\sigma}\mathcal{N}_g$ we have
$$\rho_\pi(\phi_{\infty}(a))(\xi\otimes z)=(\phi_k(a)\otimes I_{H})(\xi\otimes z)=(\phi(a)\xi)\otimes z=(\phi_k(a)\otimes I_{\mathcal{N}_g})(\xi\otimes z),$$
and
$$\rho_\pi(T_{\theta})(\xi\otimes z)=(T_{\theta}\otimes I_{H})(\xi\otimes z)=\theta\otimes\xi\otimes z=(T_{\theta}\otimes I_{\mathcal{N}_g})(\xi\otimes z).$$
Thus, the representation $\rho_\pi(H^{\infty}(E))|_{\mathcal{M}_g}$ is equal to the representation $\rho_{\sigma'}(H^{\infty}(E))=\sigma'^{\mathcal{F}(E)}(H^{\infty}(E))$.

Using the fact that $\{\hat{\pi}(M)h_{\tau}\}$ is dense in $K_{\tau}$ and $\{\sigma(M)g_0\}$ is dense in $\mathcal{N}_g$ we define for every $k\geq 0$ the operator:

$$w_k:E^{\otimes k}\otimes_{\hat{\pi}}K_{\tau} \rightarrow E^{\otimes k}\otimes_{\sigma}\mathcal{N}_g,$$
by $\xi\otimes \hat{\pi}(a)h_\tau\mapsto \xi\otimes\sigma(a)g_0$,
$\xi\in E^{\otimes k}$, $a\in M$. Since $\{\xi\otimes \hat{\pi}(a)h_\tau\}$ and $\{\xi\otimes\sigma(a)g_0\}$  span $E^{\otimes k}\otimes_{\hat{\pi}}K_{\tau}$ and $E^{\otimes k}\otimes_{\sigma}\mathcal{N}_g$ respectively, the operator $w_k$ is well defined.

For $k=0$ we have already showed that $w_0$ is an isometry from $K_{\tau}$ onto $\mathcal{N}_g$ that intertwines the representations $\hat{\pi}$ and $\sigma$.

\begin{prop}\label{w_k isometry}
The operator $w_k:E^{\otimes k}\otimes_{\hat{\pi}}K_{\tau} \rightarrow E^{\otimes k}\otimes_{\sigma}\mathcal{N}_g$ is a well defined isometry that intertwines the representation $\hat{\sigma}(\cdot)|_{E^{\otimes k}\otimes_{\hat{\pi}} K_\tau}$
and $\sigma(\cdot)|_{E^{\otimes k}\otimes_\pi H}$
\end{prop}

\noindent{\textbf{Proof.}} Let $\xi_i\otimes\hat{\pi}(a_i)h_\tau$, $i=1,2$, be in $E^{\otimes k}\otimes_{\hat{\pi}}K_{\tau}$, and
$w_k(\xi_i\otimes\hat{\pi}(a_i)h_\tau)=\xi_i\otimes\sigma(a_i)g_0$.
Denoting $c=\langle \xi_1,\xi_2\rangle$ we obtain
$$\langle \xi_1\otimes\hat{\pi}(a_1)h_\tau,\xi_2\otimes\hat{\pi}(a_2)h_\tau\rangle=
\langle\hat{\pi}(a_2)^*\hat{\pi}(c)^*\hat{\pi}(a_1)h_\tau,h_\tau\rangle=
\langle\hat{\pi}(a_2^*c^*a_1)h_\tau,h_\tau\rangle.$$

Similarly,
$$\langle \xi_1\otimes\sigma(a_1)g_0,\xi_2\otimes \sigma(a_2)g_0\rangle=
\langle\sigma(a_2^*c^*a_1)g_0,g_0\rangle,$$
so, $w_k$ is an isometry.

Let $\xi\otimes k\in E^{\otimes k}\otimes_{\hat{\pi}}K_\tau$. Then $w_k(\xi\otimes k)=\xi\otimes z\in E^{\otimes k}\otimes_{\sigma}\mathcal{N}_g$,
and
$$w_k((\phi_k(a)\otimes I_{K_\tau})(\xi\otimes k))=w_k((\phi(a)\xi)\otimes k)=(\phi(a)\xi)\otimes z.$$
But
$$(\phi(a)\xi)\otimes z=(\phi_k(a)\otimes I_{\mathcal{N}_g})(\xi\otimes z)=
(\phi_k(a)\otimes I_{\mathcal{N}_g})w_k(\xi\otimes k).$$
This proves the intertwining
$$w_k((\phi_k(a)\otimes I_{K_\tau})(\xi\otimes k))=(\phi_k(a)\otimes I_{\mathcal{N}_g})w_k(\xi\otimes k).$$
\begin{flushright}
$\square$
\end{flushright}
From the definition of the generalized powers we see that each $w_k$ is associated with $w_0$ by the identity $V^{\otimes k}(\xi)w_0=w_k\hat{V}^{\otimes k}(\xi)$, $\xi\in E^{\otimes k}$.

%\vskip 0.5 cm

Now we set
\begin{equation}\label{Operator W}
W=\sum_k w_k:\mathcal{F}(E)\otimes_{\hat{\pi}}K_\tau\rightarrow \mathcal{F}(E)\otimes_\pi H.
\end{equation}

It follows from the Propositions $\ref{Properties of w_0}$ and $\ref{w_k isometry}$ that $W$ is a well defined isometry and its image is $\mathcal{M}_g$.

\begin{rem} It is obvious from the definition of $w_k$ that $w_k(E^{\otimes k}\otimes_{\hat{\pi}}K_{\tau})=E^{\otimes k}\otimes_{\hat{\sigma}}w_0(K_{\tau})$
Fix $x\in \mathcal{F}(E)\otimes_{\hat{\pi}}K_{\tau}$ of the form $x=\xi\otimes k$,
$\xi\in \mathcal{F}(E)$, and $k\in K_\tau$.
Then we can write $Wx=W(\xi\otimes k)=\xi\otimes w_0k$. Hence, $W=I_{\mathcal{F}(E)}\otimes w_0$.
\end{rem}

\begin{prop}\label{Poperties of W}
The operator $W$ is an isometry from $\tilde{K}:=
\mathcal{F}(E)\otimes_{\hat{\pi}}K_\tau\subset G^{(\infty)}$ into $\mathcal{F}(E)\otimes_\pi H$ with $\mathcal{M}_g$ as a final subspace.
Further, $W$ intertwines the representations $\hat{\rho}|_{\tilde{K}}$ and $\rho|_{\mathcal{M}_g}$ of the algebra $H^{\infty}(E)$:
\begin{equation}\label{intertwining for W}
W\hat{\rho}(X)|_{\tilde{K}}=\rho(X) W.
\end{equation}
for every $X\in H^{\infty}(E)$.
\end{prop}

\noindent{\textbf{Proof.}}
Its remains to show only the intertwining property. To show it, it is enough to show that (\ref{intertwining for W}) holds for the generators $\{T_{\xi}, \phi_{\infty}(a):\xi\in E,a\in M\}$ of the Hardy algebra.

Since $W|_{E^{\otimes k}\otimes_{\hat{\pi}}K_\tau}=w_k$, the equality
$W\hat{\rho}(\phi_{\infty}(a))=
\rho(\phi_{\infty}(a))W$, $a\in M$, follows form Proposition $\ref{w_k isometry}$.

Now let $X=T_{\xi}$. Then $\hat{\rho}(T_{\xi})=
T_{\xi}\otimes I_{H^{(\infty})}$ and $\rho(T_{\xi})=
T_{\xi}\otimes I_H$. Fix $\eta\otimes k\in \mathcal{F}(E)\otimes_{\hat{\pi}}K_{\tau}$,
then using the previous remark we obtain\\
$W(T_{\xi}\otimes I_{H^{(\infty})})(\eta\otimes k)=
W(\xi\otimes\eta\otimes k)=
\xi\otimes\eta\otimes w_0k=
(T_{\xi}\otimes I_{\mathcal{N}_g})(\eta\otimes w_0k)=
(T_{\xi}\otimes I_H)W(\eta\otimes k)$. 
\begin{flushright}
$\square$
\end{flushright}
We obtained an isometry $W:\tilde{K}=\mathcal{F}(E)\otimes_{\hat{\pi}}K_\tau\rightarrow \mathcal{F}(E)\otimes_{\pi}H$ with final subspace $\mathcal{M}_g=\mathcal{F}(E)\otimes_{\sigma}\mathcal{N}_g$ that intertwines the induced representations $\hat{\rho}$ and $\rho$ of Hardy algebra $H^{\infty}(E)$.
In the paper ~\cite{MuS1}, partial isometries that lies in $\pi^{\mathcal{F}(E)}(\mathcal{T}_+(E))'$ are called inner operators. In our case the isometry
$W$ acts between different spaces, but intertwining $\hat{\rho}$ and $\rho$. So, it is natural to call such operators \emph{inner} operators.

We present here the general definition
\begin{defn}\label{Def of inner operator and outer vector}
Given two normal representations $\pi$ and $\mu$ of $M$ on Hilbert spaces $H$ and $K$ respectively. \\
1) An isometry $W:\mathcal{F}(E)\otimes_{\mu}K\rightarrow \mathcal{F}(E)\otimes_{\pi}H$ will be called an inner operator if

    (a) $K\subseteq H^{(\infty)}$ is a $\hat{\pi}(M)$- invariant subspace of $H^{(\infty)}$, where $\hat{\pi}$ be the ampliation of $\pi$ on $H^{(\infty)}$ and $\mu=\hat{\pi}|_{K}$. In other words, $K$ is an $M$-submodule of $H^{(\infty)}$ with respect to $\hat{\pi}$.

    (b) $W\rho_{\mu}(X)=\rho_{\pi}(X)W,\,\,\,\, X\in H^{\infty}(E)$. \\
2) A vector $y\in \mathcal{F}(E)\otimes_{\mu}K$ will be called outer if $\overline{\rho_{\mu}(H^{\infty}(E))y}=\mathcal{F}(E)\otimes_{\mu}K$.
\end{defn}
This definition and Proposition $\ref{Poperties of W}$ gives us the following Beurling type theorem for cyclic subspaces $\mathcal{M}_g$ that are considered as $\rho(H^{\infty}(E))$-modules.

\begin{thm}\label{Beurling for cyclic case}
Let $g\in G_0=\mathcal{F}(E)\otimes_\pi H$ and let
$\mathcal{M}_g=\overline{\rho((H^{\infty}(E))g}$ be a cyclic $\rho_\pi(H^{\infty}(E))$-submodule in $G_0$. Then there is a subspace $\mathcal{K}\subseteq H^{(\infty)}$ which is $\hat{\pi}(M)$-invariant and an inner operator $W:\mathcal{F}(E)\otimes_\mu \mathcal{K}\rightarrow G_0$ (where we write $\mu$ for $\hat{\pi}|_{\mathcal{K}}$) such that

1)
\begin{equation}\label{Equation Beurling for cyclic case}
\mathcal{M}_g=W(\mathcal{F}(E)\otimes_{\mu}\mathcal{K}).
\end{equation}

2) The vector $y:=W^*g$ is outer in $\mathcal{F}(E)\otimes_{\mu}\mathcal{K}$.
\end{thm}

The outer vector $y=W^*g$ will be called the outer part of $g$. Thus, the outer part of an arbitrary $g\in G$ is an outer vector in the sense of Definition $\ref{Def of inner operator and outer vector}$.

\begin{defn}\label{Definition of i-o factor of vector g} In the notation of the previous theorem, the equality
\begin{equation}\label{I-O factr of element}
g=Wy,
\end{equation}
will be called the inner-outer factorization of $g\in G_0$.
\end{defn}

The first part of the following theorem was already proved:
\begin{thm} \label{Inner-outer for element of Fock sp}
Let $\pi:M\rightarrow B(H)$ be a faithful normal representation of $W^*$-algebra $M$ on Hilbert space $H$. If $g\in\mathcal{F}(E)\otimes_{\pi}H$ then there is a $M$-submodule $\mathcal{K}\subset H^{(\infty)}$ with respect to the infinite ampliation $\hat{\pi}$ of $\pi$, an inner operator $W:\mathcal{F}(E)\otimes_{\mu}\mathcal{K}\rightarrow \mathcal{F}(E)\otimes_{\pi}H$, where $\mu=\hat{\pi}|_{\mathcal{K}}$, and an outer vector $y\in \mathcal{F}(E)\otimes_{\mu}\mathcal{K}$ such that $g=Wy$ is an inner-outer factorization of $g$.

This factorization is unique in the following sense. Let $i=1,2$ and let $\mathcal{K}_i$ are two $M$-submodules in $H^{(\infty)}$ with respect to $\hat{\pi}$ and let $\mu_i=\hat{\pi}|_{\mathcal{K}_i}$ be two normal representations of $M$ on $\mathcal{K}_i$. Suppose further that $W_i:\mathcal{F}(E)\otimes_{\mu_i}\mathcal{K}_i\rightarrow \mathcal{F}(E)\otimes_{\pi}H$ are inner operators and $y_i\in \mathcal{F}(E)\otimes_{\mu_i}\mathcal{K}_i$ are outer vectors such that $g=W_1y_1=W_2y_2$. Then there is a unitary $U:\mathcal{F}(E)\otimes_{\mu_1}\mathcal{K}_1\rightarrow \mathcal{F}(E)\otimes_{\mu_2}\mathcal{K}_2$ such that $Uy_1=y_2$ and the equality $U\rho_{\mu_1}(X)=\rho_{\mu_2}(X)U$ holds for every $X\in H^{\infty}(E)$.
\end{thm}

\noindent{Proof.} It remains to prove the uniqueness part. Let
$$W_i:\mathcal{F}(E)\otimes_{\mu_i}\mathcal{K}_i\rightarrow \mathcal{F}(E)\otimes_{\pi}H,$$
where $\mu_i$, $\mathcal{K}_i$, $y_i$, $i=1,2$, are as in the statement of the theorem.
Then $W_iy_i=g$ and
\begin{equation}\label{W_i intertwines rho_tau_i and rho_pi}
W_i\rho_{\mu_i}(X)=\rho_{\pi}(X)W_i,\,\,\ X\in H^{\infty}(E),\,\,i=1,2.
\end{equation}
Since $y_i=W^*_ig$ are outer in  $\mathcal{F}(E)\otimes_{\mu_i}\mathcal{K}_i$, $i=1,2$, and since $W_i$ have a common final subspace $\mathcal{M}_g\subset \mathcal{F}(E)\otimes_\pi H$, we get
$$W_1\overline{\mu^{\mathcal{F}(E)}_1(H^{\infty}(E))y_1}=\mathcal{M}_g=
W_2\overline{\mu^{\mathcal{F}(E)}_2(H^{\infty}(E))y_2}.$$
Set $U:=W^*_2W_1:\mathcal{F}(E)\otimes_{\mu_1}\mathcal{K}_1\rightarrow \mathcal{F}(E)\otimes_{\mu_2}\mathcal{K}_2$. Then $U$ is a unitary operator and $Uy_1=y_2$.
Finally, from the intertwining relation $(\ref{W_i intertwines rho_tau_i and rho_pi})$ we obtain
$$W_1\rho_{\mu_1}(X)W_1^*=W_2\rho_{\mu_2}(X)W_2^*.$$
Hence,
$$U\rho_{\mu_1}(X)=\rho_{\mu_2}(X)U,$$
as we wanted. $\square$

\begin{rem}
Note that, in fact, the unitary $U$ appearing in the proof can be thought of as a partial isometry in $\rho_{\hat{\pi}}(H^{\infty}(E))'$.
\end{rem}

\subsection{Inner-Outer factorization of elements of the algebra $\rho_\pi(H^{\infty}(E))'$}

We shall now apply Theorem $\ref{Inner-outer for element of Fock sp}$ to get an inner-outer factorization of an element of the commutant $\rho_\pi(H^{\infty}(E))'$.

First we consider the simple case when $\pi$ is a cyclic representation of the algebra $M$, i.e. we assume that
there is $h\in H$ such that $\overline{\pi(M)h}=H$.

Fix $S\in \rho_\pi(H^{\infty}(E))'$ and set $g:=S(1\otimes h)\in \mathcal{F}(E)\otimes_{\pi}H$, where $1\otimes h\in M\otimes_\pi H$ and $h$ is a $\pi$-cyclic vector in $H$.

Now form the subspace
$\mathcal{M}_g=\overline{\rho_{\pi}(H^{\infty}(E))g}=\overline{\rho_{\pi}(H^{\infty}(E))S(1\otimes h)}$. Since $S$ is in the commutant of $\rho_\pi(H^{\infty}(E))$ and $h$ is $\pi$-cyclic we obtain
$$\overline{\rho_{\pi}(H^{\infty}(E))S(1\otimes h)}= \overline{S\rho_{\pi}(H^{\infty}(E))(1\otimes h)}=\overline{S(\mathcal{F}(E)\otimes_\pi H)}.$$
Thus,
$$\mathcal{M}_g=\overline{S(\mathcal{F}(E)\otimes_\pi H)}.$$

By Theorem $\ref{Inner-outer for element of Fock sp}$ there are a $\hat{\pi}$-invariant Hilbert subspace $\mathcal{K}\subseteq H^{(\infty)}$, an outer element $y\in \mathcal{F}(E)\otimes_\tau\mathcal{K}$, with $\tau=\hat{\pi}|_{\mathcal{K}}$, and an inner operator
$W:\mathcal{F}(E)\otimes_{\tau}\mathcal{K}\rightarrow \mathcal{F}(E)\otimes_{\pi}H$ such that $Wy=g$ and $\mathcal{M}_g$ is the final subspace of $W$.

We set
\begin{equation}\label{Outer part of operator in commutant}
Y:=W^*S:\mathcal{F}(E)\otimes_{\pi}H\rightarrow\mathcal{F}(E)\otimes_{\tau}\mathcal{K}.
\end{equation}

\begin{prop}\label{Properties of outer op. Y in cyclic case for commutant}

1) $\overline{Y(\mathcal{F}(E)\otimes_{\pi}H)}= \mathcal{F}(E)\otimes_{\tau}\mathcal{K}$;

2) $Y\rho_\pi(X)=\rho_\tau(X)Y$, $\forall{X}\in H^{\infty}(E)$.
\end{prop}
\noindent{Proof.} 1) Since $S(\mathcal{F}(E)\otimes_\pi H)$ is dense in $\mathcal{M}_g$ and since $W$ is an isometry with $\mathcal{M}_g$ as its final subspace, we obtain that
$W^*S(\mathcal{F}(E)\otimes_\pi H)$ is dense in $\mathcal{F}(E)\otimes_\tau \mathcal{K}$.

2) Since $\rho_\pi(X)W=W\rho_\tau(X)$ and $S$ is in commutant of $\rho_\pi(H^{\infty}(E))$, we have
$$Y\rho_\pi(X)=W^*S\rho_\pi(X)=W^*\rho_\pi(X)S=\rho_\tau(X)W^*S=\rho_\tau(X)Y,$$
$X\in H^{\infty}(E)$. $\square$

The operator $Y$ will be called the outer part of $S$ and the equality $S=WY$ we call the inner-outer factorization of the operator $S\in \rho_\pi(H^{\infty}(E))'$.
The definition of the operator $Y$ a priory depends on the choice of the cyclic vector $h$. Let $h'\in H$ be another cyclic vector, $\overline{\pi(M)h'}=H$, and set $g':=S(a\otimes h')$ and $\mathcal{M}_{g'}=\overline{\rho_{\pi}(H^{\infty}(E))S(1\otimes h')}$. Then $$\mathcal{M}_{g'}=\overline{S\rho_{\pi}(H^{\infty}(E))(1\otimes h')}=\mathcal{M}_{g}.$$
Now, by Theorem $\ref{Inner-outer for element of Fock sp}$, there are $\hat{\pi}$-invariant Hilbert subspace $\mathcal{K}'\subseteq H^{(\infty)}$, representation $\tau'=\hat{\pi}|_{\mathcal{K}'}$, the outer vector $y'\in \mathcal{F}(E)\otimes_{\tau'}\mathcal{K}'$ and an inner operator $W'$ such that $W'y'=g'$. Then the corresponding outer part is $Y'=W'^*S$. The operators $W$ and $W'$ have a common final subspace $\mathcal{M}_g$ and we define $U:=W^*W'$. Hence, the operator $U:\mathcal{F}(E)\otimes_{\tau'}\mathcal{K}'\rightarrow\mathcal{F}(E)\otimes_\tau\mathcal{K}$
is unitary such that $W'=WU$ and $U\rho_{\tau'}(X)=\rho_{\tau}(X)U$. The last intertwining relation follows as in the proof of Theorem $\ref{Inner-outer for element of Fock sp}$. Further, we have $Y=W^*S$, $Y'=W'^*S$ and then $S=WY=W'Y'=WUY'$. Thus, $Y=UY'$.
This shows that the definition of $Y$ does not depend on the choice of the cyclic element $h\in H$ up to the unitary operator $U$.

Any operator $Z:\mathcal{F}(E)\otimes_{\pi}H
\rightarrow\mathcal{F}(E)\otimes_\tau\mathcal{K}$ with dense range that intertwines the representations $\rho_{\tau}$ and $\rho_\pi$ of $H^{\infty}(E)$, will be called an \emph{outer} operator. Before we give the general definition we consider the general case of noncyclic representation $\pi$.

%%%%%%%%%%%%%%%%%%General case%%%%%%%%%%%%%%%%%%%%%%%%%%%%%%%%%%%%%

So let $\pi:M\rightarrow B(H)$ be, as usual, a faithful normal representation and let $S\in \rho_\pi(H^{\infty}(E))'$

Set $\mathcal{M}:=\overline{S(\mathcal{F}(E)\otimes_{\pi}H)}$ and let
$P_{\mathcal{N}}:=P_{\mathcal{M}}-L(P_{\mathcal{M}})$ be a wandering projection with range $\mathcal{N}$. Then in terms of the shift $\mathfrak{L}$ we get
the Wold decomposition $\mathcal{M}=\mathcal{N}\oplus \mathfrak{L}(\mathcal{N})\oplus \mathfrak{L}^2(\mathcal{N})\oplus\cdots$ that we can identify with
\begin{equation}\label{Wold for N}
\mathcal{M}=\mathcal{N}\oplus (E\otimes_\sigma\mathcal{N})\oplus (E^{\otimes 2}\otimes_\sigma\mathcal{N})\oplus\cdots
\end{equation}
Consider the restriction of $\sigma(M)|_{\mathcal{N}}$. Then $\mathcal{N}$
can be written as a direct sum $\sum_i^{\oplus} \mathcal{N}_i$ of $\sigma(M)|_{\mathcal{N}}$-cyclic subspaces $\mathcal{N}_i$ with cyclic vectors $g_i\in \mathcal{N}$, such that $\mathcal{N}_i=\overline{\sigma(M)g_i}$. Thus,
$$\mathcal{N}=\sum_{i}^{\oplus}\overline{\sigma(M)g_i}.$$

The representation $(V,\sigma)$ is an isometric representation and the generalized powers
$\tilde{V}_k:E^{\otimes k}\otimes_{\sigma}(\mathcal{F}(E)\otimes_{\pi}H)\rightarrow \mathcal{F}(E)\otimes_{\pi}H$ are isometric operators. It follows that if either $k\neq l$ or $i\neq j$ one has $E^{\otimes k}\otimes_\sigma(\sigma(M)g_i)\perp E^{\otimes l}\otimes_\sigma(\sigma(M)g_j)$.

Then the Wold decomposition ($\ref{Wold for N}$) can be written as
$$\mathcal{M}=\sum^{\oplus}_i\overline{\sigma(M)g_i}\oplus
(E\otimes_\sigma\sum^{\oplus}_i\overline{\sigma(M)g_i})\oplus ...\oplus
(E^{\otimes k}\otimes_\sigma\sum^{\oplus}_i\overline{\sigma(M)g_i})\oplus ...$$

Rearranging terms we can write
$$\mathcal{M}=
\mathcal{M}_{g_1}\oplus\mathcal{M}_{g_2}\oplus...\oplus\mathcal{M}_{g_m}\oplus ...,$$
where $\mathcal{M}_{g_i}=\sum_k^{\oplus}E^{\otimes k}\otimes_\sigma\mathcal{N}_i=\overline{\rho_\pi(H^{\infty}(E))g_{i}}$ with $\mathcal{N}_i$ as a wandering subspace in $\mathcal{M}_{g_i}$ (and thus the cyclic vectors $g_i$ are wandering). From now on we shall write
$\mathcal{M}_i=\mathcal{M}_{g_i}$ and then
$\mathcal{M}=\sum_i^{\oplus}\mathcal{M}_i$.

Since all $\mathcal{M}_i$ are pairwise orthogonal we may apply Theorem $\ref{Inner-outer for element of Fock sp}$ for every $i$. So, for every $i$ there is a $\hat{\pi}(M)$-invariant Hilbert subspace $\mathcal{K}_i\subseteq H^{(\infty)}$, a normal representation $\tau_i=\hat{\pi}|_{\mathcal{K}_i}$ of $M$ on $\mathcal{K}_i$, an outer element $y_i\in \mathcal{F}(E)\otimes_{\tau_i}\mathcal{K}_i$ and an inner operator $W_i:\mathcal{F}(E)\otimes_{\tau_i}\mathcal{K}_i\rightarrow \mathcal{F}(E)\otimes_\pi H$ such that $g_i=W_iy_i$, the final subspace of $W_i$ is $\mathcal{M}_i$ and
$W_i\rho_{\tau_i}(X)=\rho_\pi(X)W_i$, $X\in H^{\infty}(E)$. Further, for every $i$ we set $Y_i:=W^*_iS$. The representation $\sigma(M)|_{\mathcal{N}}$ is cyclic when restricted to $\mathcal{N}_i$, hence by Theorem $\ref{Properties of outer op. Y in cyclic case for commutant}$ the operator $Y_i$ has a dense range in $\mathcal{F}(E)\otimes_{\tau_i}\mathcal{K}_i$, intertwines the representations $\rho_{\pi}$ and $\rho_{\tau_i}$ of $H^{\infty}(E)$, and it does not depend on the choice of the cyclic element up to some unitary operator $U_i$.
%Thus, $Y_i$ is the outer part of the restriction $S|_{\mathcal{M}_i}$.

\begin{rem}\label{Cardinality of ampliation}
Each $\mathcal{K}_i$ is a $\hat{\pi}(M)$-invariant subspace of $H^{(\infty)}$. If we write $n$ for the cardinality of the set of the cyclic vectors $\{g_i\}$, then, since $H$ is separable, $n\leq \aleph_0$. Thus, identifying $H^{(\infty)}$ with $(H^{(\infty)})^{(n)}$ we can, and will, assume that $\{\mathcal{K}_i\}$ are pairwise orthogonal subspaces in $H^{(\infty)}$ and we write $\mathcal{K}=\sum^{\oplus}_i\mathcal{K}_i$. In this case the representation $\tau=\sum_i\tau_i$ is subrepresentation of $\hat{\pi}$ obtained by restricting $\hat{\pi}$ to the $\hat{\pi}(M)$-invariant subspace $\mathcal{K}\subseteq H^{(\infty)}$
\end{rem}

In view of this remark, the operator $W:=\sum_iW_i$ acts from the subspace $\mathcal{F}(E)\otimes_\tau\mathcal{K}\subseteq \mathcal{F}(E)\otimes_{\hat{\pi}}H^{(\infty)}$ into $\mathcal{F}(E)\otimes_\pi H$ and is an inner operator. We also write $Y:=\sum_iY_i:\mathcal{F}(E)\otimes_\pi H\rightarrow \mathcal{F}(E)\otimes_\tau\mathcal{K}$, and it follows that $S=WY$.
%%%%%%%%%%%%%%%%%%%%%%%%%%%%%%%%%%%%%%%%%%%%%%%%%%%%%%%%%%%%%%%%%%%%%%%%%%%%%%%%%
%%%%%%%%%%%%%%%%%%%%%%%%%%%%%%%%%%%%%%%%%%%%%%%%%%%%%%%%%%%%%%%%%%%%%%%%%%%%%%%%%

\begin{defn}\label{Def of an outer op.} In the above notations, each operator $Y:\mathcal{F}(E)\otimes_{\pi} H\rightarrow \mathcal{F}(E)\otimes_{\tau}\mathcal{K}$ that has a dense range and such that
$Y\rho_\pi(X)=\rho_\tau(X)Y$, for every $X\in H^{\infty}(E)$, will be called an outer operator.

If $S\in\rho_\pi(H^{\infty}(E))'$ then every factorization of $S$ is of the form
\begin{equation}\label{defn of I-O of Sin commutant}
S=WY,
\end{equation}
where $Y$ is an outer operator with a dense range in $\mathcal{F}(E)\otimes_\tau\mathcal{K}$, and $W$ is an inner operator from $\mathcal{F}(E)\otimes_\tau\mathcal{K}$ into $\mathcal{F}(E)\otimes_\pi H$ will be called an inner-outer factorization of $S$. The operator $Y$ in such factorization will be called the outer part of $S$. We write also $Y_S$ for $Y$.
\end{defn}
The outer part $Y_S=W^*S$ of $S\in \rho_\pi(H^{\infty}(E))'$ is indeed an outer operator since
$$\rho_{\tau}(X)Y_S=\rho_{\tau}(X)W^*S=W^*\rho_{\pi}(X)S=W^*S\rho_{\pi}(X)=Y_S\rho_{\pi}(X).$$

We proved the existence part of the following theorem
\begin{thm} \label{Inner-Outer for commutant}
Let $S\in \rho_\pi(H^{\infty}(E))'$. Then there exist a $\hat{\pi}$-invariant subspace $\mathcal{K}\subseteq H^{(\infty)}$, a normal representation $\tau=\hat{\pi}|_{\mathcal{K}}$ of $M$ on $\mathcal{K}$, an inner
operator $W:\mathcal{F}(E)\otimes_{\tau}\mathcal{K}\rightarrow \mathcal{F}(E)\otimes_{\pi}H$ and an outer operator
$Y:\mathcal{F}(E)\otimes_{\pi}H\rightarrow \mathcal{F}(E)\otimes_{\tau}\mathcal{K}$
such that $S=WY$.

This factorization is unique in the following sense. If there is other
$\hat{\pi}$-invariant subspace $\mathcal{K}'\subseteq H^{(\infty)}$, a normal representation $\tau'=\hat{\pi}|_{\mathcal{K}'}$ of $M$ on $\mathcal{K}'$, and if
 $S=W'Y'$, where $W':\mathcal{F}(E)\otimes_{\tau'}\mathcal{K}'\rightarrow \mathcal{F}(E)\otimes_{\pi}H$
is an inner operator with final subspace $\mathcal{M}=\overline{S(\mathcal{F}(E)\otimes_\pi H)}$, and $Y':\mathcal{F}(E)\otimes_{\pi}H\rightarrow \mathcal{F}(E)\otimes_{\tau'}\mathcal{K}'$, is an outer operator,
then there exist a unitary operator $U:\mathcal{F}(E)\otimes_{\tau}\mathcal{K}\rightarrow
\mathcal{F}(E)\otimes_{\tau'}\mathcal{K}'$ such that $W'=UW$ and $Y'=U^*Y$, and
$U\rho_{\tau}(X)=\rho_{\tau'}(X)U$, $X\in H^{\infty}(E)$.
\end{thm}

\noindent{Proof.} The existence is proved above. For the uniqueness set $U=W^*W'$. Since $W$ and $W'$ have a common final subspace, the operator $U$ is unitary and $W'=UW$. From $W'Y'=S=WY$ we easily obtain that $Y'=U^*Y$. The intertwining property for $U$ follows form the definition of $U$ and from the intertwining properties of $W$ and $W'$. As in inner-outer factorization of vector, the unitary $U$ can be thought of as a partial isometry in $\rho_{\hat{\pi}}(H^{\infty}(E))'$. $\square$
\vskip 0.5cm
Let $V\in \rho_\pi(H^{\infty}(E))'$ be a partial isometry and let $V=WY$ be its inner-outer factorization. In this case the outer part of $Y$ is also a partial isometry with $\ker{Y}=\ker{V}$.

In the paper \cite{MuS1} Muhly and Solel proved Beurling Theorem for $\mathcal{T}_+(E)$-invariant subspaces. They considered the $C^*$-correspondence $E$ and assumed that $\mathcal{T}_+(E)$ is represented by some isometric representation. In their proof they used an additional assumption of quasi-invariance of the representation $\pi$, \cite[page 868]{MuS1}. J. Meyer in his Ph.D. Thesis \cite{Me} pointed out that if $\pi$ is a faithful normal representation of a $W^*$-algebra $M$, $E$ is a $W^*$- correspondence over $M$ and $\rho$ is the induced representation $\rho_\pi$ of $H^{\infty}(E)$, then the quasi-invariance assumption is fulfilled. Hence, the theorem can be formulated as follows:

\begin{thm}\label{Beurling-Muhly-Solel_I} For every $\rho(H^{\infty}(E))$-invariant subspace $\mathcal{M}$ there exist a family of partial isometries $\{V_i\}_i\subset \rho(H^{\infty}(E))'$ such that ranges of $V_i$ are pairwise orthogonal and $\mathcal{M}=\sum_i V_i(\mathcal{F}(E)\otimes_\pi H)$.
\end{thm}

Since $H$ assumed to be separable, the family $(V_i)_i$ is at most countable. Now we apply Theorem $\ref{Inner-Outer for commutant}$ for each $V_i$ to obtain an inner-outer decomposition $V_i=W_iY_i$, where $W_i:\mathcal{F}(E)\otimes_{\tau_i}\mathcal{K}_i\rightarrow \mathcal{F}(E)\otimes_{\pi}H$ is the inner operator corresponding to $V_i$.
Set as above $\mathcal{K}=\sum_i^{\oplus}\mathcal{K}_i$ and $\tau=\sum_i^{\oplus}\tau_i$. Then $\mathcal{F}(E)\otimes_{\tau}\mathcal{K}=
\sum_{\iota}\mathcal{F}(E)\otimes_{\tau_i}\mathcal{K}_i$ and write $W=\sum_iW_i$. Then the Beurling theorem of Muhly and Solel can be reformulated in the following way.
\begin{thm}\label{Beurling-Muhly-Solel_II}
Let $\pi:M\rightarrow B(H)$ be a faithful normal representation and let $\rho_\pi:X\mapsto X\otimes I_H$ be the representation induced by $\pi$ of the Hardy algebra $H^{\infty}(E)$. Further, let $\mathcal{M}\subseteq \mathcal{F}(E)\otimes_\pi H$ be a $\rho_\pi(H^{\infty}(E))$-invariant subspace. Then there exists a sequence of inner operators $W_i:\mathcal{F}(E)\otimes_{\tau_i}\mathcal{K}_i\rightarrow \mathcal{F}(E)\otimes_\pi H$ with pairwise orthogonal ranges $\{\mathcal{M}_i\}$ such that
\begin{equation}
\mathcal{M}=W(\mathcal{F}(E)\otimes_{\tau}\mathcal{K}),
\end{equation}
where $\mathcal{F}(E)\otimes_{\tau}\mathcal{K}=\sum_i^{\oplus}\mathcal{F}(E)\otimes_{\tau_i}\mathcal{K}_i$
and $W=\sum_iW_i$.
\end{thm}
\begin{rem} 1) The initial projections $V^*_i V_i$ also lie in the commutant $\rho_\pi(H^{\infty}(E))'$. Since $V_i=W_i Y_i$, then these projections are $Y^*_i Y_i$.

2) Every $\rho_\pi(H^{\infty}(E))$-invariant subspace in $\mathcal{F}(E)\otimes_\pi H$ is a direct sum of a cyclic subspaces $\mathcal{M}_{g_i}$ for some $g_i\in \mathcal{F}(E)\otimes_\pi H$, $i\in \mathbb{N}$.

\end{rem}

\subsection{Factorization of elements of  $\rho_\pi(H^{\infty}(E))$}

In this subsection we use the concept of duality of $W^*$-correspondences to produce a natural factorization of an arbitrary element of $\rho_\pi(H^\infty(E))$. This concept was developed in ~\cite[Section 3]{MuS3}.

Let $\pi: M\rightarrow B(H)$ be a normal representation of $M$ on a Hilbert space $H$. We put
\begin{equation}
E^{\pi}:=\{\eta:H\rightarrow E\otimes_{\pi}H:\eta \pi(a)=(\phi(a)\otimes
I_H)\eta, a\in M\}.
\end{equation}

On the set $E^{\pi}$ we define the structure of a $W^*$-correspondence over the von Neumann algebra $\pi(M)'$ putting
$\langle \eta,\zeta \rangle:=\eta^*\zeta$ for the $\pi(M)'$-valued inner product, $\eta,\zeta\in E^{\pi}$. It is easy to check that $\langle \eta,\zeta \rangle\in \pi(M)'$. For the bimodule operations: $b\cdot \eta=(I\otimes b)\eta$, and
$\eta\cdot c=\eta c$, where $b,c \in \pi(M)'$.
\begin{defn} The $W^*$-correspondence $E^{\pi}$ is called the $\pi$-dual of $E$.
\end{defn}

Let $\iota:\pi(M)'\rightarrow B(H)$ be the identity representation.
Then we can form $E^{\pi,\iota}:=(E^{\pi})^{\iota}$.
So, $E^{\pi,\iota}=\{S:H\rightarrow E^{\pi}\otimes_{\iota}H:
S\iota(b)=\iota^{E^{\pi}}\circ\phi_{E^{\pi}}(b)S,
b\in \pi(M)'\}$. This is a $W^*$-correspondence over $\pi(M)''=\pi(M)$.

In \cite{MuS3} it was proved that for every faithful normal representation $\pi$ of a $W^*$-algebra $M$, every $W^*$-correspondence $E$ over $M$ is isomorphic to
$E^{\pi,\iota}$.
We give a short description of this isomorphism.

For $\xi\in E$ let $L_{\xi}:H\rightarrow E\otimes_\pi H$ be defined by $L_\xi =\xi\otimes h$, $h \in H$. Then
$L_{\xi}$ is a bounded linear map since $\|L_{\xi}h\|^2\leq \|\xi\|^2\|h\|^2$
and $L_{\xi}^*(\zeta\otimes h)=\pi(\langle \xi,\zeta \rangle)h$.
For each $\xi\in E$ we define the map $\hat{\xi}:H\rightarrow E^{\pi}\otimes_{\iota}H$ by
means of its adjoint:
$$\hat{\xi}^*(\eta\otimes h)=L_{\xi}^*(\eta(h)),$$
$\eta\otimes h\in E^{\pi}\otimes_{\iota} H$.

\begin{thm}(\cite[Theorem 3.6]{MuS3})\label{Isomorphism with the second dual corresp}
If the representation $\pi$ of $M$ on $H$ is faithful, then the map $\xi\mapsto \hat{\xi}$ just defined, is an isomorphism of the $W^*$-correspondences $E$ and $E^{\pi,\iota}$.
\end{thm}

For every $k\geq 0$, let $U_k:E^{\otimes k}\otimes_{\pi}H\rightarrow (E^{\pi})^{\otimes k}\otimes_{\iota}H$
be the map defined in terms of its adjoint by
$U_k^*(\eta_1\otimes...\otimes \eta_n\otimes h)=(I_{E^{\otimes k-1}}\otimes \eta_1)...(I_{E}\otimes
\eta_{k-1})\eta_k(h)$.
It is proved in \cite{MuS3} that $U_k$
is a Hilbert space isomorphism from $E^{\otimes k}\otimes_{\pi}H$ onto $(E^{\pi})^{\otimes k}\otimes_{\iota}H$.

By Theorem $\ref{Isomorphism with the second dual corresp}$, for every $k\geq 1$ the $W^*$-correspondence $E^{\otimes k}$ over $M$ is isomorphic to the $W^*$-correspondence $(E^{\otimes k})^{\pi,\iota}\cong (E^{\pi,\iota})^{\otimes k}$. If $\xi\in E^{\otimes k}$ then the corresponding element $\widehat{\xi}\in (E^{\otimes k})^{\pi,\iota}$
is defined now by the formula
$$\widehat{\xi}^*(\eta_1\otimes...\otimes\eta_k\otimes h)=L_{\xi}^*U^*_k(\eta_1\otimes...\otimes\eta_k\otimes h),$$
where $L_{\xi}:h\mapsto \xi\otimes h$ is a bounded linear map from $H$ to $E^{\otimes k}\otimes_{\pi}H$.
Thus, we obtain
\begin{equation}\label{Formula for hat(xi)}
\hat{\xi}=U_kL_{\xi},\,\,\ \text{for}\, \xi\in E^{\otimes k}.
\end{equation}

For the dual correspondence ($\pi$-dual to $E$) we can form the (dual) Fock space $\mathcal{F}(E^{\pi})$, which is a $W^*$-correspondence over $\pi(M)'$, and the Hilbert space $\mathcal{F}(E^{\pi})\otimes_\iota H$. Let us define
$U:=\sum^{\oplus}_{k\geq 0}U_k$. It follows that
the map $U:=\sum^{\oplus}_{k\geq 0}U_k$ is a Hilbert space isomorphism from
$\mathcal{F}(E)\otimes_{\pi}H $ onto $\mathcal{F}(E^{\pi})\otimes_{\iota}H$, and its adjoint acts on decomposable tensors by
$U^*(\eta_1\otimes...\otimes\eta_n\otimes h)=(I_{E^{\otimes n-1}}\otimes \eta_1)...(I_{E}\otimes
\eta_{n-1})\eta_nh$.

\begin{defn}\label{Fourier trans definition}
The map $U_{\pi}=U:\mathcal{F}(E)\otimes_{\pi}H\rightarrow\mathcal{F}(E^{\pi})\otimes_{\iota}H$
will be called the Fourier transform determined by $\pi$.
\end{defn}

%%%%%%%%%%%%%%%%%%%%%%%%%%%%%%%%%%%%%%%%%%%%%%%%%%%%%%%%%%%%%%%%%%%%%%%%%%%%%
%%%%%%%%%%%%%%%%%%%%%%%%%%%%%%%%%%%%%%%%%%%%%%%%%%%%%%%%%%%%%%%%%%%%%%%%%%%%%%%%
Let $\pi:M\rightarrow B(H)$ be a faithful normal representation. Then there exists a natural isometric representation of $(E^{\pi},\pi(M)')$ on $\mathcal{F}(E)\otimes_{\pi}H$ induced
by $\pi$. Let $\nu:\pi(M)'\rightarrow B(\mathcal{F}(E)\otimes_{\pi}H)$
be a $*$-representation defined by $\nu(b)=I_{\mathcal{F}(E)}\otimes b$. Then $\nu$ is a faithful normal
representation of the von Neumann algebra $\pi(M)'$ and by Theorem $\ref{Rieffel thm}$,
$\pi^{\mathcal{F}(E)}(\mathcal{L}(\mathcal{F}(E)))'=\nu(\pi(M)')=\{I_{\mathcal{F}(E)}\otimes b :b\in
\pi(M)'\}$. Given $\eta\in E^{\pi}$, for each $n\geq 0$ the operators
$L_{\eta,n}:E^{\otimes n}\otimes_{\pi}H\rightarrow E^{\otimes n+1}\otimes_{\pi}H$ are defined
by $L_{\eta,n}(\xi\otimes h)=\xi\otimes \eta h$, where we have identified $E^{\otimes
n+1}\otimes_{\pi}H$ with $E^{\otimes n}\otimes_{\pi^{E}\circ \phi}(E \otimes_{\pi}H)$.
Since $\|L_{\eta,n}\|\leq \|\eta\|$, we may define the operator
$\Psi(\eta):\mathcal{F}(E)\otimes_{\pi}H\rightarrow \mathcal{F}(E)\otimes_{\pi}H$
by $\Psi(\eta)=\sum^{\oplus}_{k\geq 0}L_{\eta,k}$. Thus we may think of $\Psi(\eta)$ as $I_{\mathcal{F}(E)}\otimes \eta$ on
$\mathcal{F}(E)\otimes_{\pi}H$. It is easy to see that $\Psi$ is a bimodule map. For the
inner product, let $\eta_1, \eta_2\in E^{\pi}$ and $\xi\otimes h, \zeta\otimes k\in E^{\otimes
n}\otimes_{\pi}H$, then a simple calculation shows that

$$\langle \Psi(\eta_1)(\xi\otimes h),\Psi(\eta_2)(\zeta\otimes k)\rangle=
\langle \xi\otimes h,\nu(\eta_1^*\eta_2)(\zeta\otimes k)\rangle,$$
so, $(\Psi,\nu)$ is an isometric representation of $(E^{\pi},\pi(M)')$ on the Hilbert space
$\mathcal{F}(E)\otimes_{\pi}H$.
Combining the integrated form $\nu \times \Psi$ of $(\Psi,\nu)$ with the definition of the Fourier transform $U=U_\pi$ we obtain

\begin{equation}\label{Psi(eta)formula}
U^*\iota^{\mathcal{F}(E^{\pi})}(T_{\eta})U=\Psi(\eta),
\end{equation}
where $\eta\in E^{\pi}$ and $T_{\eta}$ is the corresponding creation operator
in $H^{\infty}(E^{\pi})$, and
\begin{equation}\label{pi_Psi(a)formula}
U^*\iota^{\mathcal{F}(E^{\pi})}(\phi_{E^{\pi},\infty}(b))U=\nu(b),
\end{equation}
where $b\in \pi(M)'$ and $\phi_{E^{\pi},\infty}$ is the left action of $\pi(M)'$ on $\mathcal{F}(E^{\pi})$.  This equality can be rewritten as
\begin{equation}\label{Intertwining M' and U}
U(I_{\mathcal{F}(E)}\otimes b)=(\phi_{E^{\pi},\infty}(b)\otimes I_H)U.
\end{equation}
Thus, the Fourier transform $U=U_\pi$ intertwines the actions of $\pi(M)'$ on $\mathcal{F}(E)\otimes_\pi H$ and on $\mathcal{F}(E^{\pi})\otimes_\iota H$ respectively.

The following theorem identifies the commutant of the Hardy algebra  represented by an induced representation.

\begin{thm}(\cite{MuS3}, Theorem 3.9)\label{commutant of ind. rep.}
Let $E$ be a $W^*$-correspondence over $M$, and let $\pi:M\rightarrow B(H)$ be a faithful
normal representation of $M$ on a Hilbert space $H$.
Write $\rho_\pi$ for the representation $\pi^{\mathcal{F}(E)}$ of $H^{\infty}(E)$ on $\mathcal{F}(E)\otimes_{\pi}H$ induced by $\pi$, and write
$\rho^\pi$ for the representation of $H^{\infty}(E^{\pi})$ on
$\mathcal{F}(E)\otimes_\pi H$ defined by
\begin{equation}\label{Formula for rho^sigma}
\rho^\pi(X)=U^*\iota^{\mathcal{F}(E^{\pi})}(X)U,
\end{equation}
with $X\in H^{\infty}(E^{\pi})$.
Then $\rho^\pi$ is an
ultraweakly continuous, completely isometric representation of $H^{\infty}(E^{\pi})$ that extends the representation $\nu\times\Psi$ of $\mathcal{T}_+(E^{\pi})$, and
$\rho^\pi(H^{\infty}(E^{\pi}))$ is the commutant of $\rho_\pi(H^{\infty}(E))$, i.e. $\rho^\pi(H^{\infty}(E^{\pi}))=\rho_\pi(H^{\infty}(E))'$.
\end{thm}

\begin{cor}(\cite{MuS3}, Corollary 3.10)\label{bicommutant of ind. rep.}
In the preceding notation,
$\rho_\pi(H^{\infty}(E))''=\rho_\pi(H^{\infty}(E))$.
\end{cor}

%%%%%%%%%%%%%%%%%%%%%%%%%%%%%%%%%%%%%%%%%%%%%%%%%%%%%%%%%%%%%%%%%%%%%%%%%%%%
%%%%%%%%%%%%%%%%%%%%%%%%%%%%%%%%%%%%%%%%%%%%%%%%%%%%%%%%%%%%%%%%%%%%%%%%%%%%
%%%%%%%%%%%%%%%%%%%%%%%%%%%%%%%%%%%%%%%%%%%%%%%%%%%%%%%%%%%%%%%%%%%%%%%%%%%%%%

\vspace{3mm}
Now we turn to the factorization of elements of $\rho_\pi(H^{\infty}(E))$. It will be obtained as a corollary of Theorem $\ref{Inner-Outer for commutant}$.

Let $X\otimes I_H\in\rho_\pi(H^{\infty}(E))$ and set
$$\mathcal{M}:=\overline{(X\otimes I_H)(\mathcal{F}(E)\otimes_\pi H)}.$$
Then $\mathcal{M}$ is $\rho_\pi(H^{\infty}(E))'$-invariant. Now let $U_\pi$ be a Fourier transform defined by $\pi$. Then the subspace
$$\tilde{\mathcal{M}}:=U_\pi\mathcal{M}\subseteq\mathcal{F}(E^\pi\otimes_\iota H)$$
is $\rho_\iota(H^{\infty}(E^\pi))$ - invariant, where by $\rho_\iota$ we denote the induced representation $\iota^{\mathcal{F}(E^{\pi})}$ of $H^{\infty}(E^\pi)$ on $\mathcal{F}(E^\pi)\otimes_\iota H$.
Set $\tilde{X}=U_\pi(X\otimes I_H)U_\pi^*$. Then $\tilde{X}$ is in the commutant of $\rho_\iota(H^{\infty}(E^\pi))$ (see Theorem $\ref{commutant of ind. rep.}$) and
$$\tilde{\mathcal{M}}=\overline{U_\pi(X\otimes I_H)U_\pi^*(\mathcal{F}(E^\pi)\otimes_\iota H)}=\overline{\tilde{X}(\mathcal{F}(E^\pi)\otimes_\iota H)}.$$
Write $\hat{\iota}$ for the ampliation of $\iota$ on the space $H^{(\infty)}$.
By Theorem $\ref{Inner-Outer for commutant}$ there is a $\hat{\iota}(\pi(M)')$-invariant subspace $\mathcal{L}$ in $H^{(\infty)}$, an inner operator
$$\tilde{W}:\mathcal{F}(E^\pi)\otimes_{\hat{\tau}}\mathcal{L}\rightarrow \mathcal{F}(E^\pi)\otimes_\iota H,$$
where $\hat{\tau}=\hat{\iota}|_{\mathcal{L}}$,
with a final subspace $\tilde{\mathcal{M}}$, and an outer operator
$\tilde{Y}=\tilde{W}^*\tilde{X}$, $\overline{\tilde{Y}(\mathcal{F}(E^\pi)\otimes_\iota H)}=\mathcal{F}(E^\pi)\otimes_{\hat{\tau}} \mathcal{L}$, such that $\tilde{X}=\tilde{W}\tilde{Y}$ is the inner-outer factorization of $\tilde{X}$.

Hence, $\tilde{X}=U_\pi(X\otimes I_H)U^*_\pi=\tilde{W}\tilde{Y}$, and
\begin{equation}\label{Preliminary I-O of X tensor I_H}
X\otimes I_H=U^*_\pi\tilde{W}\tilde{Y}U_\pi.
\end{equation}

\begin{thm}\label{Thm on I-O of X tensor I_H} For every $X\in H^{\infty}(E)$ the operator $\rho_\pi(X)=X\otimes I_H$
can be factorized as
\begin{equation}\label{I-O of X tensor I_H}
X\otimes I_H=\mathcal{W}\mathcal{Y},
\end{equation}
where $\mathcal{W}$ and $\mathcal{Y}$ satisfy
 
 1) $\mathcal{W}$ is a partial isometry from $\mathcal{F}(E^\pi)\otimes_{\hat{\iota}}H^{(\infty)}$ into $\mathcal{F}(E)\otimes_\pi H$ with intertwining relation
$$\mathcal{W}\rho_{\hat{\tau}}(S)=\rho^\pi(S)\mathcal{W},\,\,\, S\in H^{\infty}(E^\pi).$$
 
 2) $\mathcal{Y}$ acts from $\mathcal{F}(E)\otimes_\pi H$ into $\mathcal{F}(E^\pi)\otimes_{\hat{\iota}}H^{(\infty)}$
and satisfies the intertwining relation
$$\mathcal{Y}\rho^\pi(S)=\rho_{\hat{\tau}}(S)\mathcal{Y},\,\,\, S\in H^{\infty}(E^\pi).$$
 
 3) the initial subspace of $\mathcal{W}$ is the closure of the range of $\mathcal{Y}$.\\
This factorization is unique up to a multiplication by unitary.
\end{thm}

\noindent{Proof.} In $(\ref{Preliminary I-O of X tensor I_H})$ set
$$\mathcal{W}=U^*_\pi\tilde{W}\,\,\,\text{and}\,\,\,\mathcal{Y}=\tilde{Y}U_\pi.$$
We have seen that $\mathcal{W}$ is a partial isometry from $\mathcal{F}(E)\otimes_{\hat\tau}\mathcal{L}$ into $\mathcal{F}(E)\otimes_\pi H$ with the final subspace $\mathcal{M}$, and that $\mathcal{Y}$ is the operator from $\mathcal{F}(E)\otimes_\pi H$ and has a closed range in $\mathcal{F}(E)\otimes_{\hat{\tau}}\mathcal{L}$.

Since $\tilde{W}$ is inner, then
$\tilde{W}\rho_{\hat{\tau}}(S)=\rho_{\iota}(S)\tilde{W}$ for every $S\in H^{\infty}(E^\pi)$. Now, $U^*_\pi\rho_\iota(S)=\rho^\pi(S)U^*_\pi$, where $\rho^\pi(S)=U_\pi(\iota^{\mathcal{F}(E^\pi)}(S)U_\pi^*$ is the representation of $H^{\infty}(E^\pi)$ on $\mathcal{F}(E)\otimes_\pi H$ defined in $(\ref{Formula for rho^sigma})$. Thus,
$$\mathcal{W}\rho_{\hat{\tau}}(S)=\rho^\pi(S)\mathcal{W}.$$

Similarly we can show that
$$\mathcal{Y}\rho^\pi(S)=\rho_\iota(S)\mathcal{Y},\,\,\,\forall{S}\in H^{\infty}(E^\pi).$$

The uniqueness up to multiplication by unitary follows from the uniqueness of the inner-outer factorization $\tilde{X}=\tilde{W}\tilde{Y}$. $\square$
\vspace{2mm}

%%%%%%%%%%%%%%%%%%%%%%%%%%%%%%%%%%%%%%%%%%%%%%%%%%%%%%%%%%%%%%%%%%%%%%%%%%%%%%%%%%%%
%%%%%%%%%%%%%%%%%%%%%%%%%%%%%%%%%%%%%%%%%%%%%%%%%%%%%%%%%%%%%%%%%%%%%%%%%%%%%%%%%%%%

Department of Mathematics, Technion, Haifa, Israel,

email: lhelmer@tx.technion.ac.il

\end{document}